\documentclass[11pt]{article}
\usepackage{anysize}\marginsize{3.5cm}{3.5cm}{2cm}{3cm}
\usepackage{amssymb}
\newtheorem{satz}{Satz}[section]

\newtheorem{lemma}[satz]{Lemma}

\newtheorem{theorem}[satz]{Theorem}

\newtheorem{thevarthm}[satz]{\varthmname}

\newenvironment{varthm*}[1]{\trivlist\item[]{\bf #1.}\it}{\endtrivlist}

\renewcommand\emptyset{\varnothing}  % from amssymb
\renewcommand\ge{\geqslant}  % from amssymb
\renewcommand\geq{\geqslant}  % from amssymb
\renewcommand\leq{\leqslant}  % from amssymb
\renewcommand\epsilon{\varepsilon}
\renewcommand\phi{\varphi}

\newcommand\be{\begin{eqnarray*}}
\newcommand\ee{\end{eqnarray*}}
\def\ben{\begin{eqnarray}}
\def\een{\end{eqnarray}}

\newcommand\eqnref[1]{(\ref{#1})}

\newcommand\eps{\varepsilon}
\renewcommand\P{\mathbb P}

\newcommand\R{\mathbb R}
\newcommand\Z{\mathbb Z}

\newcommand\lra{\longrightarrow}
\newcommand\newop[2]{\def#1{\mathop{\rm #2}\nolimits}}
\newop\upper{upper}
\newcommand\calo{{\cal O}}

\newcommand\mult{{\rm mult}}

\begin{document}

\title{Seshadri constants and surfaces of minimal degree}
\author{Wioletta Syzdek and Tomasz Szemberg}
\date{\today}
\maketitle
\thispagestyle{empty}

%*****************************************************************************

\begin{abstract}
   In \cite{SyzSze} we showed that if the multiple point Seshadri
   constants of an ample line bundle
   on a smooth projective surface in very general points satisfy
   certain inequality then the surface is fibred by curves computing
   these constants. Here we characterize the border case of polarized
   surfaces whose Seshadri constants in general points fulfill the
   equality instead of inequality
   and which are not fibred by Seshadri curves. It turns
   out that these surfaces are the projective plane and surfaces of
   minimal degree.
\end{abstract}

%*****************************************************************************

\section*{Introduction and the main result}

\begin{plaintheoremnumbers}

   Given a smooth projective variety $X$ and a nef line bundle $L$ on
   $X$, Demailly defines the Seshadri constant of $L$
   at a point $P\in X$ as the real number
   $$
      \eps(L;P):=\inf_C\frac{L.C}{\mult_PC} \ ,
   $$
   where the infimum is taken over all reduced and irreducible
   curves passing through $P$
   (see \cite{Dem92} and \cite[Chapt.~5]{PAG}).

   This concept was extended by Xu \cite{Xu95}
   to finite subsets of a given
   variety. Let $r$ be an integer and $P_1,\dots, P_r$ points in
   $X$. Then the $r$-tuple Seshadri constant of $L$ at the set
   $P_1,\dots,P_r$ is the real number
      $$
      \eps(L;P_1,\dots,P_r):=\inf_{C\cap\{P_1,\dots, P_r\}\neq\emptyset}
      \frac{L.C}{\sum\mult_{P_i}C} \ ,
   $$
   where the infimum is taken over all irreducible curves passing
   through at least one of the points $P_1,\dots,P_r$.

   There is an alternative and useful description of Seshadri
   constants in terms of the nef cone of a blown up variety.
   Specifically, let
   $f:Y\lra X$ be the blowing up of $P_1,\dots,P_r\in X$
   with exceptional divisors $E_1,\dots,E_r$.
   Then the Seshadri constant can be computed as
   $$
      \eps(L;P_1,\dots,P_r)=
      \sup\left\{\lambda>0:\, f^*L-\lambda\cdot\sum_{i=1}^rE_i \mbox{ is nef}\right\}.
   $$

   The Kleiman criterion of ampleness implies then that the multiple
   point Seshadri constants are subject to the following upper bound
   which depends only on the degree of $L$ and the number of points
   $$
      \eps(L;P_1,\dots,P_r)\leq\sqrt[\dim X]{\frac{L^{\dim X}}{r}}=:\alpha(L;r).
   $$

   Whenever there is a strong inequality
   \ben\label{submax}
      \eps(L;P_1,\dots,P_r)<\alpha(L;r)
   \een
   then the Seshadri constant is actually computed by a curve and
   not approximated by a sequence of curves. For
   Seshadri constants at a single point
   this follows from
   \cite[Lemma 5.2]{Bau99}
   and the argument easily modifies to the multiple point case.
   We call any curve $C$ with
   $$\eps(L;P_1,\dots,P_r)=\frac{L.C}{\sum\mult_{P_i}C}$$
   a {\it Seshadri curve for $L$ at the $r$-tuple $P_1,\dots,P_r$}.

   Oguiso (see \cite{Ogu02}) studied the behavior of Seshadri
   constants $\eps(L;P)$ under the variation of the point $P$. He
   showed that the Seshadri function
   $$
   \eps_1:X\ni P\lra \eps(L;P)\in\R
   $$
   is
   semi-continuous and that it attains its maximal value at a set
   which is a complement of an at most countable union of Zariski
   closed proper subsets of $X$ i.e. for a very general point $P$.
   Oguiso arguments can be easily adapted to finite subsets.
   By $\eps(L; r)$ we will
   abbreviate the maximal value of the function
   $$
   \eps_r:X^r\ni(P_1,\dots,P_r)\lra\eps(L;P_1,\dots,P_r)\in\R
   $$
   i.e. $\eps(L;r):=\max\eps_r$.

   Nakamaye (see \cite[Corollary 3]{Nak03}) observed that in case of
   surfaces an inequality of type
   $$\eps(L;1)<\lambda \cdot \alpha(L;1)$$
   with a small factor $\lambda$ has strong consequences for the
   geometry of the surface. Namely there exists a non-trivial
   fibration of $X$ over a curve $B$ whose fibers are Seshadri
   curves for $L$. On surfaces this was studied in more detail by
   Tutaj-Gasi\'nska and the second author \cite{SzeTut04}. Hwang and
   Keum passed from surfaces to varieties of arbitrary dimension
   (see \cite{HwaKeu03}). In \cite{SyzSze} we started research along
   the same lines for multiple point Seshadri constants. In
   particular we proved the following theorem.

\begin{varthm*}{Theorem on fibrations}\label{syzszeineq}
   Let $X$ be a smooth projective surface, $L$ a nef and big line
   bundle on $X$ and $r\geq 2$ a fixed integer. If
   \ben\label{star}
      \eps(L;r)<\sqrt{\frac{r-1}{r}}\cdot\alpha(L;r)
   \een
   then there exists a fibration $f:X\lra B$ over a curve $B$ such
   that given $P_1, \dots, P_r\in X$ very general,
   for arbitrary
   $i=1,\dots, r$ the fiber
   $f^{-1}(f(P_i))$ computes $\eps(L;P_1,\dots,P_r)$
   i.e. the fiber is a Seshadri curve of $L$.
\end{varthm*}

   Furthermore we showed that the bound in the Theorem is sharp in
   the sense that for every integer $r$ there exists a surface $X$
   together with an ample line bundle $L$ such that one has
   equality in \eqnref{star} and $X$ is not fibred by Seshadri curves
   of $L$.

   The purpose of this note is to characterize the pairs $(X,L)$ for
   which one has an equality in \eqnref{star} and $X$ is not fibred
   by Seshadri curves. The description of such pairs is  provided in
   the next theorem which is our main result.

\begin{theorem}\label{equality}
   Let $r\ge 2$ be a given integer, $X$ a smooth projective surface
   and $L$ a nef and big line bundle on $X$ such that
   \ben\label{rownosc}
      \eps(L;r)=\sqrt{\frac{r-1}{r}}\cdot\alpha(L;r)\;.
   \een   
   If $X$ is not fibred by Seshadri curves for $L$,
   then
   \begin{itemize}
      \item[a)] either $r=2$, $X=\P^2$ and $L=\calo(1)$,
      \item[b)] or $X$ is a surface of minimal degree in $\P^r$ and
      $L=\calo_X(1)$.
   \end{itemize}
\end{theorem}

\begin{varthm*}{Remarks}\rm

   (i) A similar theorem for $r=1$ was already obtained by us in
   \cite[Theorem 3.2]{SyzSze} but the result and the methods
   are somewhat different.

   (ii) A smooth surface is of minimal degree if and only if
   it is the Veronese
   surface in $\P^5$ or a rational normal scroll. This was proved
   by Del Pezzo (see \cite{Del1886}).
   
   (iii) The converse of the Theorem holds: for any surface $X$ of minimal
   degree and $L=\calo_X(1)$, the equality \eqnref{rownosc} holds. This is 
   easy to see taking the hyperplane section through $r$ given points.
\end{varthm*}
\end{plaintheoremnumbers}

\section{Useful Lemmas}
   Here we recall two Lemmas which are essential for the proof of
   the main result.

   The first Lemma goes back to Xu \cite[Lemma 1]{Xu95}.
\begin{lemma}\label{xu}
  Let $X$ be a smooth projective surface, let
  $(C_t,(P_1)_t,\dots,(P_r)_t)_{t\in \Delta}$ be
  a non-trivial
  one parameter family of pointed reduced and irreducible curves on $X$ and let $m_i$
  be positive integers such that $\mult_{(P_i)_t}C_t\geq m_i$ for all $i=1,\dots, r$. Then
  $$
     \begin{array}{lcl}
     \mbox{for } r=1 \mbox{ and } m_1\geq 2 && C_t^2\geq m_1(m_1-1)+1 \mbox{ and }\\
     &&\\
     \mbox{for } r\geq 2 && C_t^2\geq \sum_{i=1}^r m_i^2 - \min\{m_1,\dots,m_r\}.
     \end{array}
  $$
\end{lemma}
   The second lemma was obtained by K\"uchle in \cite{Kue96} and has purely
   arithmetical character.
\begin{lemma}\label{kuechle}
   Let $r\geq 2$ and $m_1,\dots ,m_r\in {\Z}$ be
   integers with $m_1\geq\dots\geq m_r\ge 1$ and $m_1\ge 2$. Then we
   have
   $$(r+1)\sum_{i=1}^r m_i^2 > \left(\sum_{i=1}^r m_i\right)^2+m_r(r+1).$$
\end{lemma}

\section{Proof of the Theorem}
   In this section we prove Theorem \ref{equality}. First we give a short
   overview of the proof. Since the Seshadri constants of the line bundle in
   Theorem \ref{equality} are strictly less than the upper bound, they must be
   computed by Seshadri curves.

   We investigate properties of these curves in three steps.
   First we show that
   under assumptions of Theorem \ref{equality}
   the multiplicities of Seshadri curves
   in points $P_1,\dots,P_r$ must all be equal to $1$.
   This is an arithmetical part of the proof.

   In the second step which is more analytical, we show that
   Seshadri curves must be rational.

   The third step is geometrical and realizes Seshadri curves as
   hyperplane sections of $X$ embedded in a projective space
   as a surface of minimal degree.

   Let us now turn to the details.

\subsection{Multiplicities of Seshadri curves}
   By assumptions of the Theorem \ref{equality}
   inequality \eqnref{submax} is
   satisfied so for every $r$-tuple $P_1,\dots,P_r$ there exists
   a Seshadri curve $(C;P_1,\dots,P_r)$.
   By \cite[Proposition 1.3]{Syz07} there are finitely many such curves
   for every $r$-tuple. For a very general $r$-tuple we
   have the equality
   \ben\label{sesheq}
   \frac{L.(C;P_1,\dots,P_r)}{\sum_{i=1}^r
   \mult_{P_i}C}=\frac1r\cdot\sqrt{(r-1)L^2}.
   \een
   The number of algebraic families of curves satisfying this
   equality is at most countable. So at least one of these families
   must not be discrete. From now on we are interested in Seshadri
%\fillin{is one parameter enough for argument in (C)?}
   curves $(C_t;(P_1)_t,\dots,(P_r)_t)$
   for $L$ moving in a non-trivial family over some algebraic
   set $\Delta$. Let $m_i$ be the biggest integers such that
   $$\mult_{(P_i)_t}C_t\geq m_i$$
   for all $t\in \Delta$. Making $\Delta$ a little bit smaller if necessary we
   may assume that actually $m_i=\mult_{(P_i)_t}C_t$ for all $t$.

   Renumbering the points if necessary we may also assume that
   $$m_1\geq\dots\geq m_r.$$
   There are the following three cases possible:
   \begin{itemize}
      \item[(A)] $m_r\geq 1$ and $m_1\geq 2$;
      \item[(B)] $m_1=\dots m_r=1$;
      \item[(C)] $m_r=0$.
   \end{itemize}
   In this step we want to exclude (A) and (C).

   In case (A) we are in the position to apply Lemma \ref{kuechle}.
   Thus
   $$
   \frac{1}{r+1}\left(\sum_{i=1}^r m_i\right)^2<\sum_{i=1}^r
   m_i^2-m_r\leq C_t^2,
   $$
   where the second inequality is assured by Lemma \ref{xu}.
   Multiplying the above inequality by $L^2$ and applying the index
   theorem on the right hand side we arrive to the following
   inequality
   $$
   \frac{1}{r+1}\left(\sum_{i=1}^r m_i\right)^2\cdot L^2<(L.C_t)^2.
   $$
   Dividing by the sum of multiplicities and revoking
   \eqnref{sesheq} we obtain
   $$
   \frac{1}{r+1}\cdot L^2<\frac{r-1}{r^2}\cdot L^2
   $$
   which is not possible.

   In case (C) if $r\geq 3$, then we have
   $$
   \frac{L.(C;P_1,\dots,P_{r-1})}{\sum_{i=1}^{r-1}m_i}=
   \frac{L.(C;P_1,\dots,P_r)}{\sum_{i=1}^{r}m_i}
   =\sqrt{\frac{r-1}{r}}\cdot\alpha(L;r)<\sqrt{\frac{r-2}{r-1}}\cdot\alpha(L;r-1).
   $$
   Hence our Theorem on fibrations shows that $X$ is covered by
   Seshadri curves for $L$ contradicting the assumption of Theorem
   \ref{equality}.

   If $r=2$, then by assumption we have
   $$\eps(L;1)=\sqrt{\frac14L^2}$$
   and in this case we get the same contradiction by
   \cite[Theorem]{SzeTut04}.

   Thus we showed that for $P_1,\dots,P_r$ very general the Seshadri
   curve for $L$ has multiplicities equal $1$ at all these points.
   In particular we conclude from \eqnref{sesheq} that
   \ben\label{degc2}
   L.(C;P_1,\dots,P_r)=\sqrt{(r-1)L^2}\; .
   \een
   Together with the index theorem we get
   \ben\label{degc}
   C^2\leq r-1.
   \een

\subsection{Rationality of Seshadri curves}
   In this part we follow basically the deformation argument of
   \cite {EL93} with necessary modifications. First we observe that
   one can fix the points $P_1,\dots,P_{r-1}$ and consider Seshadri
   curves for the $r$-tuples $P_1,\cdots, P_{r-1},P$ with the last
   point moving. Among these curves one can find again a non-trivial
   family $(C_t;P_1,\dots,P_{r-1},P_t)$ over some smooth base $\Delta$.
   For $t$ general the corresponding Kodaira-Spencer map
   $$T_t\Delta\lra H^0(C_t,N_{C_t/X})$$
   factorizes in fact over $H^0(C_t,N_{C_t/X}(-P_1-\dots-P_{r-1}))$.

   Lemma \ref{xu} implies that $C_t^2\geq r-1$. In view of
   \eqnref{degc} we obtain that in fact $\deg N_{C_t/X}=C_t^2=r-1$.
   Since the image of the Kodaira-Spencer map is non-zero we
   conclude that the line bundle $N_{C_t/X}(-P_1-\dots-P_{r-1})$ is
   trivial. Equivalently, there is a section $s_r$ in $H^0(C_t,N_{C_t/X})$
   whose zero locus is exactly the divisor $P_1+\dots+P_{r-1}$.
   Fixing $P_r$ and moving instead another point in the tuple we get
   in the same manner sections $s_1,s_2\dots,s_r$ in $H^0(C_t,N_{C_t/X})$
   whose zero loci are $P_2+\dots +P_r$, $P_1+P_3+\dots +P_r$,
   $\dots$, $P_1+\dots +P_{r-1}$ respectively.
   They are obviously
   independent. This shows that $N_{C_t/X}$ is a line bundle of
   degree $r-1$ with at least $r$ sections. This can happen only in
   the case when $C_t$ is a rational curve. Thus we showed that
   under assumptions of Theorem \ref{syzszeineq} the Seshadri curves
   are rational.

\subsection{Embedding $X$ as a surface of minimal degree}
   It follows from the last part that $X$ is rationally connected
   hence it is a rational surface. Since $C^2=r-1$ for Seshadri curves,
   it follows from the index theorem and \eqnref{degc2} that the
   Seshadri curves are numerically equivalent. On rational surfaces
   this implies the linear equivalence, so Seshadri curves move in a
   single linear system. We call this system $|M|$ and we
   show that $M$ is in fact very ample.

   First we show that $M$ separates points.
   Let $P$ and $Q$ be two distinct points on $X$. Let
   $C$ be a Seshadri curve for $L$ lying in $|M|$ and passing
   through $P$. It might happen that $Q$ lies also on $C$. Taking
   $P_2,\dots,P_{r-1}$ general on $C$ we have that
   $(C;P,P_2,\dots,P_{r-1},Q)$ is a Seshadri curve for $L$. Taking
   $Q'$ very general away of $C$ there exists also a Seshadri curve
   $(C';P,P_2,\dots,P_{r-1},Q')$. Since $C.C'=r-1$ this new curve
   cannot pass through $Q$ and thus we separated $P$ and $Q$.

   Next we show that $M$ separates tangent vectors. To this end for
   a fixed point $P$ it
   is enough to find two Seshadri curves intersecting transversally
   at $P$. Again, this is the case for the curves $C$ and $C'$ from
   the argument above as they have $r-1=C.C'$ points in common, so
   must intersect at every of these points transversally.

   If $r=2$ then $M$ has degree $1$. This shows that $X$ is $\P^2$.
   For $r\geq 3$ and a smooth curve $C\in |M|$ we consider the exact sequence
   $$0\lra \calo_X\lra \calo_X(M)\lra \calo_C(C)\lra 0.$$
   Since $H^1(\calo_X)=0$ and $h^0(C,\calo_C(C))=r$ as already
   established in the previous part, we conclude
   from the long cohomology sequence that $M$ has $r+1$
   sections. Hence the image of $X$ under the mapping given by $|M|$
   must be a surface of minimal degree.

\paragraph*{Acknowledgement.}
   During the conference ''Linear systems and subschemes'' held in Gent
   in April 2007 we
   presented the results of \cite{SyzSze} and began investigations
   which led to the present note. We would like to thank the
   organizers of this conference for providing a nice and
   stimulating atmosphere. We would also like to thank Brian
   Harbourn for nice discussions there.

   Both authors were partially supported by
   MNiSW grant N~N201 388834.
%*****************************************************************************

%*****************************************************************************

\bigskip
\small
   Wioletta Syzdek,
   Tomasz Szemberg

   Instytut Matematyki AP,
   ul. Podchor\c a\.zych 2,
   PL-30-084 Krak\'ow, Poland

   {\it current address:}

   Mathematisches Institut, Universit\"at
   Duisburg-Essen, 45117 Essen, Germany

\nopagebreak
   E-mail: \texttt{syzdek@ap.krakow.pl}

   E-mail: \texttt{szemberg@ap.krakow.pl}
\end{document}